\documentstyle[newlfont]{amsart}

\newtheorem{thm}{Theorem}[section]
\newtheorem{lem}[thm]{Lemma}

\newtheorem{cor}[thm]{Corollary}

\theoremstyle{definition}
\newtheorem{defn}{Definition}

\newtheorem{ex}{Example}

\theoremstyle{remark}

\numberwithin{equation}{section}
\def\Ker{\operatorname {Ker}}
\def\Im{\operatorname {Im}}

\def\Aut{\operatorname {Aut}}

\def\Comod{\operatorname {Comod}}

\newcommand{\C}{\Bbb{C}}

\newcommand{\Q}{\Bbb{Q}}

\newcommand{\at}{\tilde{a}}
\newcommand{\bt}{\tilde{b}}
\newcommand{\ct}{\tilde{c}}
\newcommand{\dt}{\tilde{d}}

\newcommand{\npm}{\frak{n}^{\pm}}
\newcommand{\bpm}{\frak{b}^{\pm}}

\newcommand{\g}{\frak{g}}

\newcommand{\hgoth}{\frak{h}}
\newcommand{\h}{\frak{h}}

\newcommand{\pair}[2]{\langle {#1},{#2} \rangle}

\newcommand{\bP}{\bold{P}}

\newcommand{\bQ}{\bold{Q}}
\newcommand{\bB}{\bold{B}}
\newcommand{\bR}{\bold{R}}
\newcommand{\bRp}{\bold{R}_+}

\newcommand{\bL}{\bold{L}}

\newcommand{\calCq}{\cal{C}_q}
\newcommand{\uqbm}{U_q(\frak{b}^-)}
\newcommand{\uqbp}{U_q(\frak{b}^+)}

\newcommand{\qhat}{(q-q^{-1})}

\newcommand{\uqg}{U_q(\g)}
\newcommand{\uqdg}{U_q(\frak{d}(\g))}
\newcommand{\cqg}{{\Bbb{C}}_q[G]}
\newcommand{\cqdg}{{\Bbb{C}}_q[D(G)]}
\newcommand{\cqds}{{\Bbb{C}}_q[D(SL(2))]}

\newcommand{\op}{{\text{op}}}
\newcommand{\cop}{{\text{cop}}}

\newcommand{\pr}[2]{\langle #1, #2 \rangle}

\begin{document}


\title{Double quantum groups and Iwasawa decomposition}
\author{Timothy J. Hodges} 
\thanks{The author was partially supported by grants from the
National 
Security Agency and the National Science Foundation.}
\address{University of Cincinnati\\ Cincinnati, OH 45221-0025\\
U.S.A.}
\email{timothy.hodges@@uc.edu}


\begin{abstract}
     The double quantum groups $\cqdg = \cqg \Join \cqg$ are the
Hopf algebras underlying the complex quantum groups of which the
simplest example is the quantum Lorentz group. They are non-
standard quantizations of the double group $G \times G$. We
construct a corresponding quantized universal enveloping algebra
$\uqdg$ and prove that the pairing between $\cqdg$ and $\uqdg$ is
nondegenerate. We analyze the representation theory of these
$\cqdg$, give a detailed version of the Iwasawa decomposition
proved by Podles and Woronowicz for the quantum Lorentz group,
and show that $\cqdg$ is noetherian. Finally we outline how to
construct more general non-standard quantum groups using quantum
double groups and their generalizations.
\end{abstract}

\maketitle


\section{Introduction}

     A {\em double quantum group} is a Hopf algebra of the form
$A
\Join A$ where $A$ is a braided Hopf algebra and $A\Join A$
denotes the generalized Drinfeld double constructed using the
pairing defined by the braiding. This construction originated
independently in Podles and Woronowicz's construction of the
quantum Lorentz group and in early work of Majid on
bicrossproducts. Here we analyze this construction from an
algebraic point of view with emphasis on the case when $A$ is the
standard quantum group $\cqg$ for $G$ a connected semi-simple
complex algebraic group. In this case the double $\cqdg =
\cqg \Join \cqg$ can be thought of as a non-standard
quantization of $\C[G \times G]$. One reason why the double is a
particularly natural object, is the existence of surjective Hopf 
algebra maps $m : \cqdg \to \cqg$ and $\theta: \cqdg \to \uqg$.
These maps are quantizations of the natural embeddings of $G$ and
its dual $G_r$ into the double group $G \times G$ that occur in
the analysis of the symplectic leaves of $G$ and of dressing
transformations on $G$ \cite{STS,LW,HL1}.

     The main results we prove here are:
\begin{enumerate}
\item The category of right $\cqdg$-comodules is equivalent to
the
category of right $\C_q[G\times G]$-comodules as a braided rigid
monoidal category.
\item The natural pairing between $\cqdg$ and its
Fadeev-Reshetikhin-Takhtadjan (FRT) dual $\uqdg$ is 
non-degenerate. Hence $\uqdg$ has a category of finite
dimensional
left
modules which is equivalent to the category of right 
$\cqdg$-comodules and $\cqdg$ is the restricted Hopf dual of
$\uqdg$ with
respect to this category.
\item The map $(m \otimes \theta)\Delta : \cqdg \to \cqg \otimes
\uqg^{\cop}$ is injective. This is a generalization of the
``Iwasawa decomposition'' proved for the quantum Lorentz group in
\cite{PW}. More precisely we show that there is a finite group
$\Gamma$ acting on $\cqg \otimes \uqg^{\cop}$ and an Ore set
$\cal{S}$ of $\cqdg$ such that the induced map from
$\cqdg_{\cal{S}}$ to $(\cqg \otimes \uqg^{\cop})^\Gamma$ is an
isomorphism.
\item The algebra $\cqdg$ is Noetherian.
\end{enumerate}
A key starting point is the observation by Doi and Takeuchi
\cite{DT} that the double $A \Join A$ can be realized as a twist
of the usual tensor product by a 2-cocycle (a construction dual
to Drinfeld's gauge transformation \cite{D1}). From this result,
assertion (1) follows easily and then (2) and the injectivity of
the Iwasawa decomposition map follow by standard calculations.
The second part of (3) (the identification of the image of the
Iwasawa decomposition map) appears to be new even in the case $G
= SL(2)$. The proof that $\cqdg$ is noetherian uses the approach
of Brown and Goodearl \cite{BG}. Since both $\cqg$ and $\uqg$ are
homomorphic images of $\cqdg$ this provides an interesting
unified proof of the noetherianity of the standard quantum
function algebra and the standard quantized enveloping algebra.

     The algebra $\C_q[SL(2)]$ was introduced and studied (using
different constructions) by Podles and Woronowicz in \cite{PW}
and
by Carow-Watamura, Schlieker, Scholl and Watamura in \cite{CSSW}.
Further work on this and more general ``complex quantum groups''
was done by a variety of authors, (see for instance \cite{CW,Ta,
DSWZ}). In \cite{M1}, Majid pointed out that these algebras could
be
constructed using the quantum double construction that we
consider
here. This approach was also taken up in \cite{CEJS}. This
article
is, in a sense, a continuation of the approach of \cite{M1,CEJS}
with
the important distinction that we do not consider the *-algebra
structure. All of the above papers are concerned with quantizing
the algebra of complex valued coordinate functions on a simple
complex algebraic group $G$, considered as a real Lie group. Thus
the relevant notion is a complex Hopf algebra equipped with a 
*-operation. Here we consider these algebras rather as
quantifications of the algebra of regular functions on the double
complex algebraic group $G \times G$. However, the underlying
Hopf
algebra is the same and  many of the basic structural results are
relevant in both contexts.

     We also show how to apply the same twisting technique to
construct some new non-standard quantum groups. Recall that the
Poisson group structures on a semi-simple Lie group $G$ given by
a
solution of the modified classical Yang Baxter equation, have
been
classified
by Belavin and Drinfeld \cite{BD}. They are given (in part) by
triples of the form $(\tau, \bB_1,\bB_2)$ where $\bB_i$ are
subsets
of the base $\bB$ of the root system for $\g$ and $\tau$ is a
bijection from $\bB_1$ to $\bB_2$ satisfying certain conditions.
In
the special case where $\bB_1$ and $\bB_2$ are ``disjoint'' we
construct quantizations of the associated Poisson group.
This construction is closely related to the construction of
non-standard quantum groups by Fronsdal and Galindo using the
universal $T$-matrix \cite{FG}.

     The Hopf algebra notation that we use is generally the
standard notation used in Sweedler's book except that we remove
the parenthesis from the `Sweedler notation' writing $\Delta(x) =
\sum x_1 \otimes x_2$ rather than $\sum x_{(1)} \otimes x_{(2)}$.
We use
$m,\Delta, \epsilon$ and $S$ indiscriminately for the
multiplication,
comultiplication, counit and antipode in any Hopf algebra. 

     This work originated in discussions between the author and
T. Levasseur on
non-standard quantum groups while the former was visiting the
Universit\'{e} de Poitiers. He would like to thank Levasseur for
his
significant contribution and the Mathematics Department in
Poitiers for their hospitality.

\section{Braided Hopf algebras and their doubles}

     An appropriately general setting in which to view the
basic constructions is that of a {\em braided Hopf algebra}.
Braided Hopf algebras are also known as dual quasi-triangular,
coquasi-triangular or cobraided Hopf algebras. They play the role
(roughly) of the dual of a quasi-triangular Hopf algebra and
include all of the standard, multi-parameter and non-standard
quantizations of semi-simple algebraic groups.  The notation
follows closely that of Doi and Takeuchi but many of the ideas
occur in the work of Majid \cite{M1} and Larson and Towber
\cite{LT}.

     Let $A$ and $U$ be a pair of Hopf algebras. A {\em skew
pairing} on the ordered pair $(A,U)$ is a bilinear form $\tau$
satisfying 
\begin{enumerate}
     \item $\tau(bc, u) = \sum\tau(b,u_1)\tau(c,u_2)$
     \item $\tau(b,uv) = \sum\tau(b_1,v)\tau(b_2,u)$
     \item $\tau(1,u) = \epsilon(u)$; $\tau(a,1)=\epsilon(a).$
\end{enumerate}
A bilinear form on the pair $(A,U)$ is said to be invertible if
it is invertible
as an element of $(A \otimes U)^*$. Any  skew pairing $\tau$
is invertible with inverse given by
$$\tau^{-1}(a,u) = \tau(S(a), u).$$
For any bilinear form $\tau$ we denote
by $\tau_{21}$ the pairing given by $\tau_{21}(a,b) =
\tau(b,a)$. Note that if $\tau$ is a skew pairing, then so is
$\tau^{-1}_{21}$.

     Given a skew pairing on a pair of Hopf algebras
$(A,U)$, one may construct a new Hopf algebra, the {\em quantum
double} $A \Join_\tau U$ (or just $A \Join U$). As a coalgebra
the double is
isomorphic to $A \otimes U$. The multiplication is given by the
rule:
$$(1 \otimes u)(a \otimes 1) = \sum \tau(a_1,u_1)a_2 \otimes u_2
\tau^{-1}(a_3,u_3)$$
or, equivalently,
$$\sum (1\otimes u_1)(a_1 \otimes 1) \tau(a_2,u_2) = \sum
\tau(a_1,u_1) a_2 \otimes u_2$$

 A {\em braided Hopf algebra} is a Hopf algebra $A$ together
with an invertible skew pairing $\beta$ on $(A,A)$ such that
$$\sum b_1a_1 \beta(a_2,b_2) = \sum \beta(a_1,b_1) a_2b_2.$$
The pairing $\beta$ is called a {\em
braiding} on $A$. Note that the antipode in a braided algebra is
always bijective \cite{Do}. Given a braided  Hopf algebra $A$ we
may form the double Hopf
algebra $A\Join_\beta A$. 

     A {\em 2-cocycle} on a Hopf algebra $A$ is an invertible
pairing $\sigma : A \otimes A \to k$ such that for all $x$, $y$
and $z$ in $A$,
$$\sum \sigma(x_{1},y_{1}) \sigma(x_{2}y_{2},z) = \sum
\sigma(y_{1},z_{1})
     \sigma(x,y_{2}z_{2}) $$
and $\sigma(1,1) = 1$.
      Given a 2-cocycle $\sigma$ on a Hopf algebra $A$, one can
twist the multiplication to get a new Hopf algebra $A_\sigma$
\cite{DT}.
The new multiplication is given by
$$x \cdot y = \sum \sigma(x_{1},y_{1}) x_{2} y_{2}
\sigma^{-1}(x_{3},y_{3}). $$
The comultiplication in $A_\sigma$ remains the same. In
particular, $A$ and $ A_\sigma$ are isomorphic as coalgebras.  
This construction is essentially the dual of Drinfeld's gauge
transformation \cite{D1}.

     Recall that the category of right comodules over a braided
Hopf algebra has a natural structure of a braided rigid monoidal
category \cite{LT}. A 2-cocycle $\sigma$ on a braided Hopf
algebra can be
used to define an equivalence of such categories between $\Comod
A$ and $\Comod A_\sigma$ when $A_\sigma$ is equipped with the
braiding defined below.

\begin{thm}
     Let $A$ be a braided  Hopf algebra with braiding $\beta$.
Let $\sigma$ be a 2-cocycle on $A$. Let $A_\sigma$  be the
twisted Hopf algebra defined above. Then $\sigma_{21}* \beta*
\sigma^{-1}$ (convolution product) is a braiding on
$A_\sigma$.
Moreover the categories of comodules over $A$ and $A_\sigma$ are
equivalent as rigid braided monoidal categories.
\end{thm} 

\begin{pf} See \cite{M2}.
\end{pf}

     Doi and Takeuchi observed in \cite{DT} that the quantum
double may be constructed from the tensor product by twisting by
a 2-cocycle.
Let $\tau$ be an invertible skew pairing between the Hopf
algebras $B$ and $H$.
 Then the  bilinear form $[\tau]$ defined on $B\otimes H$ by
$$[\tau](b\otimes g, c\otimes h) =
\epsilon(b)\epsilon(h)\tau(c,g)$$
is a 2-cocycle and the twisted Hopf algebra $(B \otimes
H)_{[\tau]}$ is
isomorphic as a Hopf algebra to the quantum double
$B\Join_{\tau} H$ \cite[Prop. 2.2]{DT}. In particular, if $(A,
\beta)$ is a braided
Hopf algebra,
then the  double $A\Join _\beta A$ is isomorphic to $(A
\otimes A)_{[\beta]}$.

     There are a number of different ways  of defining a braiding
on $A \Join A$. Recall that if $\beta$ is a braiding on $A$, then
so is $\beta_{21}^{-1}$. The tensor algebra $A \otimes A$ can
therefore be made into a braided algebra in a number of different
ways using these two braidings. It turns out that the appropriate
choice is  the braiding given by $\beta$ on the first component
and
$\beta_{21}^{-1}$ on the second component which we shall denote
by
$(\beta,\beta_{21}^{-1})$. Using the above theorem
we deduce that 
$$\gamma = [\beta]_{21} * (\beta,\beta_{21}^{-1})*[\beta]^{-1}$$
is a braiding on $A \Join A$ (cf. \cite[Proposition 2]{M1}). This
braiding makes the category of right
comodules into a braided category.

\begin{thm} \label{comod}
Let $(A, \beta)$ be a braided Hopf algebra. Then the double
$A\Join
A$ has a braiding $\gamma$ given by
$$\gamma(a \otimes b, a' \otimes b')  =
\beta(a,a_1'b_1')\beta^{-1}(a_2'b_2',b). $$
The category
$\Comod A\Join A$ is equivalent as a braided rigid monoidal
category to $\Comod A\otimes A$.
\end{thm}

     There are two important Hopf algebra homomorphisms
associated to
the double $A \Join A$ of a braided Hopf algebra $A$. The first
is the multiplication map  \cite[Prop. 3.1]{DT}
$$m:A\Join A \to A.$$
The second is the map  $$\theta : A \Join A \to
(A^\circ)^{\cop}$$ defined
by  \cite[Thm. 3.2]{DT}
$$\theta(x\otimes y ) = \beta(x,-)\beta^{-1}(-,y).$$
Notice that $\theta= m(l^+\otimes l^-)$ where $l^\pm :
A \to
(A^\circ)^{\cop}$ are the maps $l^+(x) = \beta(x,-)$ and $l^-(y)
= \beta^{-1}(-,y)$.
We denote the image of $\theta$ in $A^\circ$ by $U(A)$ and refer
to it as the
{\em FRT-dual} of $A$. Recall that $U(A)$ is said to be dense in
$A^\circ$ if the
pairing between $A$ and $U(A)$ is non-degenerate. Combining the
above maps via
the comultiplication yields an algebra map
$$ \xi = (\theta \otimes m)\Delta : A \Join A \to
U(A)^{\text{cop}}\otimes A.$$

 Given the braiding $\gamma$ on $A \Join A$ we have an associated
FRT-dual $U(A
\Join A)$ and maps $l^{\pm}: A \Join A \to U(A \Join A)$. 
For any
Hopf algebra $B$ we shall denote by $\pr{-}{-}$
the pairing between $B^{\circ}$ and $B$. 
The maps $m$ and $\theta$ have duals $m^*: U(A) \to (A
\Join A)^\circ$ and
$\theta^* : A^{\op} \to (A \Join A)^\circ$ given by
$\pr{m^*(u)}{b\otimes b'} = \pr{u}{m(b\otimes b')}= \pr{u}{bb'}$
and $\pr{\theta^*(a)}{b \otimes b'} = \pr{\theta(b\otimes
b')}{a}= \sum \beta(b,a_1)\beta^{-1}(a_2,b')$.  

\begin{lem}
Consider the maps $l^{\pm}: A \Join A \to U(A\Join A)$. Then $l^+
= m^* \theta$ and $l^- = S \theta^* m$.
In particular, the images of both $m^*$ and $\theta^*$ are
contained in $U(A \Join A)$.
\end{lem}

\begin{pf} Observe that
$$
\gamma(a \otimes b, a' \otimes b') = \pr{\theta(a\otimes b)}{m(a'
\otimes b')}
 = \pr{m^* \theta(a\otimes b)}{a' \otimes b'}.$$
Similarly,
$$
\gamma^{-1}(a \otimes b, a' \otimes b') = \pr{\theta(a\otimes
b)}{S^{-1}m(a' \otimes b')}
= \pr{S\theta^* m(a' \otimes b')}{a\otimes b} $$
\end{pf}

     Recall that the usual pairing between $U(A)$ and $A$ becomes
a skew pairing between $U(A)$ and $A^{\op}$, allowing us to form
the double $U(A)\Join A^{\op}$.

\begin{thm}\label{iwaua}
     The map 
$$\zeta = m(m^* \otimes \theta^*): U(A)\Join 
A^{\op} \to U(A
\Join A)$$
 is a  surjective Hopf algebra map. 
\end{thm}

\begin{pf} The lemma implies that the image
of $\zeta$ is precisely $U(A\Join A)$. The fact that $\zeta$ is a
Hopf algebra map is more or less well-known (see for instance
\cite{M1}). It follows from the formula
$$
 \sum \theta^*(a_1)m^*(u_1) \pr{u_2}{a_2} = \sum \pr{u_1}{a_1}
m^*(u_2) \theta^*(a_2)
$$
which may be verified directly.
\end{pf}

\begin{cor} \label{inj1}
     If $U(A\Join A)$ is dense in $(A\Join A)^\circ$, then the
map 
$$ \xi = (m \otimes \theta)\Delta : A \Join A \to A \otimes
U(A)^{\cop}
$$
is injective.
\end{cor}

\begin{pf} Notice that the pairing between $A \Join A$ and
$U(A)\Join  A^{\op}$ induced from $\zeta$ is the same as that
induced from the map $\xi$. Thus
 $$  \Ker \xi \subset (U(A)\Join 
A^{\op})^\perp \subset U(A\Join A)^\perp .$$
The density of $U(A\Join A)$ implies that $U(A\Join
A)^\perp =0$, so
$  \Ker \xi =0$ also.
\end{pf}

We now give an alternative construction of the braiding on $A
\Join
A$. This approach is analogous to Drinfeld's construction of a
universal $R$-matrix for the Drinfeld double $H^* \Join H$. We
follow the reformulation of this construction given in
\cite{Ga,Jo}.
This construction is also needed  later in defining the braiding
on
the standard
quantum groups.
Let $A$ and $B$ be Hopf algebras and let $\sigma$ be a skew
pairing
between $A$ and $B$.
Define Hopf algebra maps
$$ \Phi_1: A^{\cop} \to B^\circ, \quad  \Phi_2: B^{\op} \to
A^\circ$$
by $\Phi_1(a) = \sigma(a,-)$ and $\Phi_2(b) = \sigma(-,b)$. If 
the pairing $\sigma$ is non-degenerate these maps will be
injective. Henceforth we shall assume that this is the case. Now
suppose that we have a Hopf pairing $\phi$ between $C$ and 
$A \Join B$ which identifies $C$ with a subalgebra of the dual of
$A \Join B$. There are Hopf algebra maps
$$ \theta_1 : C \to B^\circ, \quad \theta_2 : C \to A^\circ$$
defined by $\pr{\theta_1(c)}{b}= \phi(c,1\otimes b)$ and
$\pr{\theta_2(c)}{a}=\phi(c,a\otimes 1)$ respectively. 

     In order to construct a braiding on $C$ we need to
assume that $\Im \theta_i \subset \Im \Phi_i$ for $i =1$ and $2$.
In this case we can construct maps
$$ \psi_1 = \Phi_1^{-1} \theta_1 : C \to A^{\cop}, \quad \psi_2 =
\Phi_2^{-1}\theta_2 : C \to B^{\op}$$

\begin{thm} \label{gait}
The form $ \beta  \in (C \otimes C)^*$ defined by 
$$ \beta(x,y) = \sigma(\psi_1(x), \psi_2(y))$$
is a braiding on $C$.

Let $\pi: A \Join B \to C^*$ be the natural map and let
$l^\pm : C \to U(C)^{\cop}$ be the maps described above. Then
$l^+ = \pi \psi_1$ and $l^- = S \pi \psi_2$. 
\end{thm}

\begin{pf} See \cite[9.4.6]{Jo}
\end{pf}

We now apply this result in the case where $\sigma$ is the
natural
skew pairing between $U(A)$ and $A^\op$ for some braided Hopf
algebra
$A$ and $C$ is  $A \Join A$. When $U(A\Join A)$ is dense in
$(A\Join A)^\circ$, this produces a braiding on $A \Join
A$ which we will denote by $\gamma'$.

\begin{thm}
The braidings $\gamma$ and $\gamma'$ coincide.
\end{thm}

\begin{pf} 
The pairing between $U(A)\Join A^{\op}$ and $A \Join A$ is given
by
$$\pr{u \otimes c}{a \otimes b} = \pr{m^*(u)}{a_1\otimes b_1}
\pr{\theta^*(c)}{a_2\otimes b_2}
=\pr{u}{a_1b_1}\pr{\theta(a_2\otimes b_2)}{c}$$
So
$$\pr{\theta_1(a\otimes b)}{c} =\pr{1\otimes c}{a \otimes
b}=\pr{\theta(a\otimes b)}{c} 
=    \pr{ \Phi_1 \theta(a\otimes b)}{c}.$$
Hence $\psi_1 =\Phi_1^{-1}\theta_1 =\theta$. Similarly,
$$\pr{\theta_2(a\otimes b)}{u} =\pr{u\otimes 1}{a \otimes
b}=\pr{u}{m(a\otimes b)} 
=    \pr{ \Phi_2 m(a\otimes b)}{u}.$$
Thus $\psi_2 = m$. Hence
$$
\gamma'(a \otimes b, a' \otimes b') = \pr{\theta(a\otimes
b)}{m(a'
\otimes b')}
= \gamma(a \otimes b, a' \otimes b')$$
as required.
\end{pf}

     This theorem says that the braiding $\gamma$ is essentially
the same object as the universal $T$-matrix of \cite{RS} and
\cite{FG}.

There is a natural algebra embedding of $U(A \Join A)$ into $U(A)
\otimes U(A)$ given by 
$$ \chi(u) = \sum f(u_1) \otimes f'(u_2) $$
where $f(u)$, $f'(u)$ denotes the restriction of $u$ to the first
and second copies of $A$ respectively. Of particular importance
is
the composition of $\chi$ and $\zeta$.

\begin{lem} \label{chir}
\begin{enumerate}
\item For all $a \in A$, $\chi \theta^*(a) = \sum l^-S(a_1)
\otimes
l^+S(a_2)$.
\item For all $u \in U(A)$, $\chi m^*(u) = \sum u_1 \otimes u_2 =
\Delta(u)$.
\end{enumerate}
\end{lem}

\begin{pf} Notice that
\begin{align*}
\pr{\theta^*(a)}{b \otimes c} & = \pr{\theta(b \otimes c)}{a}
          = \beta(b, a_1)\beta^{-1}(a_2,c) \\
     & = \beta^{-1}(b,S(a_1))\beta(S(a_2),c)
          = \pr{l^-S(a_1)}{b}\pr{l^+S(a_2)}{c}
\end{align*}
This proves the first assertion. The second assertion is proved
similarly.
\end{pf}

\section{Double quantum groups and the quantum Lorentz group}

     We now define the double quantum group associated to a
connected complex semi-simple algebraic group $G$.  We continue
to use the notation of \cite{HLT}. For the convenience of the
reader we recall briefly the relevant details.

Let $\g$ be the Lie algebra of $G$. Let $\h$ be a Cartan
subalgebra of $\g$, $\bR$ the associated root
system, $\bB = \{ \alpha_1, \dots, \alpha_n\}$ a basis of $\bR$,
$\bRp$ the set of positive roots and $W$ the Weyl group. 
We denote by $\bP$ and $\bQ$ the
lattices of weights and roots respectively and by $\bP^+$ the set
of
dominant integral weights. Let $H$ be a maximal torus of $G$ with
Lie algebra
$\hgoth$ and denote by $\bL$ the character group of $H$, which we
shall identify with a sublattice of $\bP$ containing  $\bQ$.
Let $(-,-)$ be  the Killing form on
$\h^*$ and set $d_i = (\alpha_i,\alpha_i)/2$. 
Set $\npm = \oplus_{\alpha \in \bold{R}_+} \g_{\pm
\alpha}, \quad \bpm = \h
\oplus \npm$.

Let $q \in \C^*$ and assume that {\em $q$ is not a root of
unity.} Since we need to consider rational powers of $q$ we adopt
the following convention. Pick $\hbar \in \C$ such that $q=e^{-
\hbar /2}$ and define $q^m = e^{-m \hbar /2}$ for all $m \in \Q$.
We set $q_i = q^{d_i}$ for $i = 1, \dots, n$.
Denote by $U^0$ the group algebra of  $\bL$,
$$U^0 = \C[k_\lambda \, ; \, \lambda \in \bL], \quad k_0 =1,
\quad
k_\lambda k_\mu = k_{\lambda + \mu}.$$
Set $k_i = k_{\alpha_i}$, $ 1 \le i \le n$. The standard
quantized
enveloping algebra associated to this data  is the Hopf
algebra 
$$U_q(\g) = U^0[e_i, f_i \, ; \, 1 \le i \le n]$$
with defining  relations:
\begin{gather*}
k_\lambda e_j k_\lambda^{-1} = q^{(\lambda,\alpha_j)} e_j, \quad
k_\lambda f_j k_\lambda^{-1} = q^{-(\lambda,\alpha_j)} f_j \\
 e_if_j - f_je_i = \delta_{ij}(k_i-k_i^{-1})/(q^{d_i}-q^{-d_i}) 
\end{gather*}
and the quantum Serre relations as given in \cite{HLT}. The Hopf
algebra structure is given
by
\begin{gather*}
\Delta(k_\lambda) = k_\lambda \otimes k_\lambda, \quad
\epsilon(k_\lambda) = 1, \quad S(k_\lambda) = k_\lambda^{-1} \\
\Delta(e_i) = e_i \otimes 1 + k_i \otimes e_i, \quad 
\Delta(f_i) = f_i \otimes k_i^{-1} + 1 \otimes f_i \\
\epsilon(e_i) = \epsilon(f_i) =0, \quad S(e_i) = -k_i^{-1} e_i,
\quad 
S(f_i) = -f_i k_i.
\end{gather*}
Notice that $\uqg$ depends on $G$ rather than $\g$. Thus the
notation is a little ambiguous.

We define subalgebras of $\uqg$ as follows
\begin{gather*}
\uqbp = U^0[e_i \, ; \, 1 \le i \le n], \quad \uqbm =  U^0[f_i
\, ; \, 1 \le i \le n].
\end{gather*}
Notice that
$U^0$ and $U_q(\frak{b}^{\pm})$ are Hopf subalgebras of $\uqg.$

Let $M$ be a left $\uqg$-mo\-du\-le. An element $x
\in M$ is said to have {\em weight $\mu \in
\bL$} if $k_\lambda x = q^{(\lambda,\mu)}x$ for all $\lambda \in
\bL$;
 we denote by $M_\mu$ the subspace of elements of weight $\mu$.
It is known that the category of finite dimensional
$\uqg$-modules is a completely reducible braided rigid monoidal
category.
Set $\bL^+ = \bL \cap \bP^+$ and recall that for each $\lambda
\in
\bL^+$ there exists a finite dimensional simple module of highest
weight $\lambda$, denoted by $L(\lambda)$. One has $L(\lambda)^*
\cong L(w_0 \lambda)$ where $w_0$
is the longest element of $W$.
Let $\calCq$ be the subcategory of finite dimensional
$\uqg$-modules consisting of finite direct sums of $L(\lambda)$,
$\lambda
\in \bL^+$. The category  $\calCq$ is closed under tensor
products and the
formation of duals. 

Let $M$ be an object of $\calCq$, then $M
=\oplus_{\mu \in \bold{L}} M_\mu$. For $f \in M^*$, $v \in M$ we
define the
 coordinate function  $c_{f,v} \in \uqg^*$ by
$$\forall u \in \uqg, \quad c_{f,v}(u) = \pair{f}{uv}$$
where $\pair{\,}{\,}$ is the duality pairing. 
 The quantized function algebra $\cqg$ is the restricted
dual of $\uqg$ with respect to $\calCq$. That is,
$$\cqg = \C[c_{f,v} \, ; \, v \in M, f \in M^*, \, M \in
\text{obj}(\calCq)].$$ 
The algebra $\cqg$ is a Hopf algebra. If $\{v_1,
\dots,
v_s; f_1, \dots, f_s\}$ is a dual basis for $M$
one has
\begin{equation} \label{formulas 3.1}
\Delta(c_{f,v}) = \sum_i c_{f,v_i} \otimes c_{f_i,v}, \quad
\epsilon(c_{f,v}) = \pair{f}{v}, \quad S(c_{f,v}) = c_{v,f}.  
\end{equation}
Notice that we may assume that $v_j \in M_{\nu_j}, \, f_j \in
M^*_{-\nu_j}$. When $v \in L(\lambda)_\mu$ and $f \in
L(\lambda)^*_{-\nu}$ we denote the element $c_{f,v}$ by
$c^\lambda_{-\nu,\mu}$. Although convenient this notation is a
little ambiguous and some care has to be taken in interpreting
the standard formulas such as 
$\Delta(c^\lambda_{\nu,\mu})=c^\lambda_{\nu,\mu_i}\otimes
c^\lambda_{-\mu_i,\mu}$ and $S(c^\lambda_{\nu,\mu}) = c^{-
w_0\lambda}_{\mu,\nu}$.
 
     Recall that the Rosso form $\phi(\;,\;)$ defines a skew
pairing on $(\uqbm,\uqbp)$. Consider the induced maps  
$$\Phi_1: \uqbm^{\cop} \to \uqbp^\circ, \quad \Phi_2 : 
\uqbp^{\op} \to \uqbm^\circ $$
and the maps
$$ \theta_1: \cqg \to \uqbp^\circ, \quad  \theta_2: \cqg \to
\uqbm^\circ $$
as defined above. By \cite[Proposition 4.6]{HLT}, we have that
$\Im \theta_i = \Im \Phi_i$. Hence the Rosso form induces a
braiding on $\cqg$ defined by
$$ \beta(x,y) = \phi(\psi_1(x),\psi_2(y))$$
where $\psi_i = \Phi_i^{-1} \theta_i$.
     
\begin{defn} 
    We define $\cqdg$ to be the double quantum group $\cqg
\Join \cqg$.
\end{defn}

\begin{ex}
In the case when $G= SL(2,\C)$, the double quantum group
$\cqds$ is the Hopf algebra underlying the quantum Lorentz group.
It is often denoted by $SL_q(2,\C)$. 
As an algebra it is generated by the elements $a,b,c,d$ and $\at,
\bt,\ct,\dt$ subject to the relations:
$$ba = q ab, \; ca = q ac, \; db = q bd, \; dc = qcd $$
     $$ cb = bc, \quad ad - da = \qhat bc, \quad da - qbc = 1$$
$$\bt  \at  = q \at \bt , \; \ct \at  = q \at c, \; \dt \bt  = q
\bt \dt , \; \dt \ct  = q\ct \dt  $$
     $$ \ct \bt  = \bt \ct , \quad  \at \dt  - \dt \at  = \qhat
\bt \ct , \quad \dt \at  - q\bt \ct  = 1$$
$$a \at = \at a, \; q a \bt = \bt a, \;
a\ct = \ct a + \qhat c \at, \; a \dt + \qhat c \bt = \dt a $$
$$qb \at = \at b = \qhat \bt a, \; b \bt = \bt b, \;
b \ct + \qhat d \at = \qhat \dt a + \at b $$
$$ b \dt + \qhat d \bt = q \dt b$$
$$c \at = q \at c, \quad c \bt = \bt c, \quad c \ct = \ct c,
\quad q c \dt = \dt c$$
$$ d \at = \at d + \qhat \bt c, \; d \bt = q \bt d, \; q d \ct
=\ct d + \qhat \dt c, \; d \dt = \dt d$$
\end{ex}

     We begin by stating explicitly Theorem \ref{comod} for this
case.  
 
\begin{thm} \label{comodG}
The category
$\Comod \cqdg$ is equivalent as a braided rigid monoidal
category to $\Comod \C_q[G\times G]$.
\end{thm}

     The equivalence of the comodule categories in the case of
the quantum Lorentz group was proved by Podles and Woronowicz in
their original paper \cite{PW}. To our knowledge no more general
result has appeared in the literature, although this result is
more or less implicit in some of Majid's work.

 \begin{defn}
    Define $\uqdg$ to be the FRT-dual $U(\cqdg)$.
\end{defn}
 
Notice again that $\uqdg$ depends not only on $\g$ but also on
the choice of $G$. Theorem 
\ref{iwaua} gives a weak version of Iwasawa decomposition for
$\uqdg$.

\begin{thm} \label{iwauqg} 
     The map
$$ \zeta = m(m^* \otimes \theta^*): \uqg\Join  \cqg^{\op} \to
\uqdg$$
is an epimorphism of Hopf algebras.
\end{thm}

     We conjecture that this map $\zeta$ is in fact an
isomorphism.

     We now show that the pairing between $\cqdg$ and $\uqdg$ is
non-degenerate.  We need a detailed description of the maps
$l^\pm$. Note that Theorem \ref{gait} implies that
$$\Ker l^+ = \{ c_{f,v} \mid f \in (\uqbp v)^\perp \}$$
and
 $$\Ker l^- = \{ c_{f,v} \mid f \in (\uqbm v)^\perp \}.$$
Hence $c^\lambda_{-\mu,\nu} \in \Ker l^+$ if  $\mu-\nu$ is not a
non-negative integer combination of positive roots. Similarly
$c^\lambda_{-\mu,\nu} \in \Ker l^-$ if  $\nu-\mu$ is not a
non-negative integer combination of positive roots.

Define
$$ U^+  = \C[e_i \mid 1 \leq i \leq n], \quad U^- = \C[f_i \mid 1
\leq i \leq n]$$
and for all $\beta \in \bQ$ set 
$$U^\pm_\beta = \{ u \in U^\pm \mid k_\lambda u k _{-\lambda} =
q^{(\lambda, \beta) } u \}.$$

\begin{lem} \label{lpm}
 Let $\lambda \in \bL^+$ and let $\mu \in \bL$.
\begin{enumerate}
\item There exists a unique $y \in U^-_{\mu-\lambda}$ such that
$l^+(c^\lambda_{-\lambda,\mu}) = y k_{-\mu}$.
\item There exists a unique $x \in U^+_{\lambda-\mu}$ such that
$l^-(c^\lambda_{-\mu,\lambda}) = x k_{\mu}$
\end{enumerate}
\end{lem}

\begin{pf} See \cite{HLT} or \cite[9.2.11]{Jo}.
\end{pf}

     Consider the natural functor from right $\cqdg$-comodules to
left $\uqdg$-modules. For a right
$\cqdg$-comodule $V$, we define an action of $\uqdg$ by
$$u.v = \sum \pr{v_1}{u}v_0$$
for $v \in V$ and $u \in \uqdg$. 

     The proof of the following theorem generalizes the proof of
an
analogous result for the quantum Lorentz group given in
\cite{Ta}.

\begin{thm}
 Let $V$ be a simple right $\cqdg$-comodule and
consider
$V$ as a left $\uqdg$-module. Then $V$ is simple.
\end{thm}

\begin{pf}
     By Theorem \ref{comodG} any irreducible $\cqdg$-comodule is
of the form
$V=L(\nu) \otimes L(\nu')$ where $\nu,\nu' \in \bL^+$.  The
action
of $\uqdg$ on $V$ is then given by the map
$\chi$ and the usual action of $\uqg \otimes \uqg$. 
Suppose that $\lambda \in \bL^+$ and $\mu \in \bL$. It follows
from Lemma \ref{chir} that for any $c \in \cqg$,
$$\chi \theta^*(c) = l^-S(c_1) \otimes l^+S(c_2)$$
Hence
$$\chi \theta^*(c^\lambda_{-\mu,\lambda}) =\sum
l^-S(c^\lambda_{-\mu,\nu}) \otimes l^+S(c^\lambda_{-\nu,\lambda})
=
l^-S(c^\lambda_{-\mu,\lambda}) \otimes l^+S(c^\lambda_{-
\lambda,\lambda})$$
and similarly
$$\chi \theta^*(c^\lambda_{-\lambda, \mu}) =
l^-S(c^\lambda_{-\lambda,\lambda}) \otimes l^+S(c^\lambda_{-
\lambda, \mu})$$
(cf. \cite[9.2.14]{Jo}). In particular Lemma \ref{lpm} yields
$$\chi\theta^*(c^\lambda_{-\lambda,\lambda}) = k_{-
\lambda} \otimes k_{\lambda}$$
 Thus the image of $\chi\theta^*$ contains an `antidiagonal' copy
of the
subalgebra $\C[k_\lambda \mid \lambda \in \bL^+]$ of 
$U^0$. On the other hand, the image under $m^*$ of $U^0$ is a
`diagonal' copy of $U^0$. It follows that any $\uqdg$-submodule
$V'$ of
$V$ is a sum of its $U^0 \otimes U^0$-weight spaces. Now by
applying Lemma \ref{lpm} again for $\mu = \lambda -
\alpha_j$ yields (up to a scalar factor)
$$\chi\theta^*(c^\lambda_{-\mu,\lambda}) = e_j k_{-\lambda}
\otimes k_{\lambda}$$
$$\chi\theta^*(c^\lambda_{-\lambda,\mu}) =  k_{-
\lambda} \otimes f_j k_{\lambda}$$
Hence if  $V'$ is a non-zero submodule of $V$ it must contain
a weight vector $v_+\otimes v_-'$ where $v_-'$ is a highest
weight
vector of $V(\nu')$ and $v_+$ is a lowest weight vector
of $V(\nu)$. However
this element generates $V$ as a $\uqg$ module
via the diagonal action. Thus $V'=V$ as required.
\end{pf}

     The theorem yields a form of Peter-Weyl theorem linking
$\cqdg$ and $\uqdg$. Denote by $\tilde{\cal{C}}_q$ the full
subcategory of $\uqdg$-mod consisting of direct sums of modules
of
the form
$L(\nu) \otimes L(\nu')$. Since $\tilde{\cal{C}}_q$ is closed
under tensor products, duals and direct sums we may form the
restricted dual of $\uqdg$ with respect to $\tilde{\cal{C}}_q$.
This is the Hopf algebra of coordinate functions
$$\uqdg^\circ_{\tilde{\cal{C}}_q} = \C[c_{f,v} \, ; \, v \in M, f
\in M^*, \, M \in
\text{obj}(\tilde{\cal{C}}_q)].$$ 

\begin{cor}
The pairing between $\cqdg$ and $\uqdg$ is non-degenerate. The
categories $\cal{C}_q$ and $\tilde{\cal{C}}_q$ are equivalent and
$\cqdg$
is the restricted dual of $\uqdg$ with respect to
$\tilde{\cal{C}}_q$.
\end{cor}

\begin{pf} 
For $\nu \in \bL^+$ let $C(L(\nu)) \subset \cqg$ be the
subcoalgebra of coordinate functions on $L(\nu)$. Then
$$ \cqdg = \oplus_{\nu,\nu' \in \bL^+} C(L(\nu)) \otimes
C(L(\nu')).$$
Similarly let $C(L(\nu) \otimes L(\nu'))$ be the subcoalgebra of
$\uqdg^*$ of coordinate functions on $L(\nu) \otimes L(\nu')$.
Then 
$$\uqdg^\circ_{\tilde{\cal{C}}_q} = \oplus_{\nu,\nu' \in \bL^+}
C(L(\nu) \otimes L(\nu')).$$
Let $\eta: \cqdg \to \uqdg^*$ be the natural map. Then since
$L(\nu) \otimes L(\nu')$ is simple $\eta$ maps $C(L(\nu)) \otimes
C(L(\nu'))$ isomorphically onto $C(L(\nu) \otimes
L(\nu'))$.
Hence $\eta$ is an isomorphism from $\cqdg$ to
$\uqdg^\circ_{\tilde{\cal{C}}_q}$.
\end{pf}

     In particular this implies that the intersection of the
annihilators of the simple modules $L(\nu)\otimes L(\nu')$ is
zero.
Hence $\uqdg$ is semi-primitive and residually finite.

     Applying Corollary \ref{inj1} yields the following version
of Iwasawa decomposition for $\cqdg$.

\begin{cor}
     The map
$$\xi=(m \otimes \theta)\Delta : \cqdg \to 
\cqg \otimes \uqg^{\text{cop}}  $$
 is injective.
\end{cor}

     This result is proved for the quantum Lorentz group by
Podles and Woronowicz in \cite[Theorem 1.3]{PW}. In \cite{Ta},
Takeuchi defines the quantum Lorentz group as a subalgebra of
$\C_q[SL(2)] \otimes U_q(\frak{sl}(2))^{\text{cop}}$ and then
goes on to prove that
this subalgebra is isomorphic to $\C_q[D(SL(2))]$. A discussion
of this map in a fairly general context appears in \cite{M1}.
Although
$\xi$ is never surjective, it is in a certain sense not far away
from being surjective. We now show that the image has a
localization which is the
invariant subring for the action of a certain finite group,
$\Gamma$. This is only to be expected if we view this map as the
quantization of the map $G \times G_r \to G \times G$ described
in \cite{HL1}. This map is a finite morphism onto the open subset
$GG_r$ of
$G \times G$ whose fibre at every point is $\Gamma$. 

     Let $\Gamma = \{ h \in H \mid h^2 = e \}$. We may identify
the dual group $\hat{\Gamma}$ with the quotient $P /2P$. 

\begin{lem} $\cqdg /(\Ker m + \Ker \theta) \cong
\C[\hat{\Gamma}]$.
\end{lem}

\begin{pf}
Denote by $\C[\Gamma]$ the group algebra of $\Gamma$. Then there
is a surjective Hopf algebra map $\eta: \cqdg \to
\C[\hat{\Gamma}] = \C[\Gamma]^*$ given by
$$\eta(c^\lambda_{-\mu,\lambda}\otimes
c^{\lambda'}_{-\mu',\lambda'})(h) =
\lambda(h)\lambda(h')\epsilon(c^\lambda_{-\mu,\lambda}
c^{\lambda'}_{-\mu',\lambda'}).$$
It is easily checked that $\Ker \eta \supset \Ker \theta + \Ker
m$. Conversely, let $t_\lambda$ be the image of $c^\lambda_{-
\lambda,\lambda}\otimes 1$ in $\cqdg/(\Ker m + \Ker \theta)$.
Then
$t_\lambda^2 = 1$, $t_\lambda t_\mu = t_{\lambda+\mu}$ and the
$t_\lambda$ span $\cqdg/(\Ker m + \Ker \theta)$. Thus by
comparing dimensions we obtain that $\Ker \eta = \Ker \theta +
\Ker m$.
\end{pf}

     Denote by $R(A)$ the group of one-dimensional
representations of a Hopf algebra $A$. The group $\Gamma$ embeds
in
$R(\cqg)$, $R(\uqg^{\cop})$ and $R(\cqdg)$ in such a way that the
embeddings commute with the induced maps $R(\uqg^{\cop}) \to
R(\cqdg)$ and $R(\cqg) \to R(\cqdg)$. Recall that for any Hopf
algebra $A$, there are left and right translation maps $l,r :
R(A) \to \Aut(A)$ given by 
$$l_h(x) = \sum h^{-1}(x_1)x_2,\quad r_h(x) = \sum x_1h(x_2).
$$  
The left and right translation actions of $\Gamma$ on
$\cqg$ and $\uqg^{\cop}$ factor through left and right
translation
actions on $\cqdg$. Define an action $\tilde{d}: \Gamma \to \Aut
\cqdg \otimes \cqdg$ by
$$ \tilde{d}_h (x \otimes y) = r_h(x) \otimes
l_h(y).$$
Consider $\Gamma$ acting similarly on $(\cqg\otimes
\uqg^{\cop})$. 

\begin{lem}    The image of the comultiplication $\Delta: \cqdg
\to
\cqdg
\otimes \cqdg$ is contained in the subring of invariants, $(\cqdg
\otimes \cqdg)^\Gamma$. Hence the image of $\cqdg$ under $\xi$ is
contained in the invariants $(\cqg\otimes \uqg^{\cop})^\Gamma$.
\end{lem}

\begin{pf} Notice that for $c_{\mu,\lambda} \in \cqg$ and $h \in
\Gamma$,
$$ l_h(c_{\mu,\lambda}) = \mu(h) c_{\mu,\lambda}, \quad
     r_h(c_{\mu,\lambda}) = \lambda(h) c_{\mu,\lambda}.$$
Hence for $c_{\mu,\lambda} \otimes c_{\mu',\lambda'} \in \cqdg$,
\begin{align*}
\tilde{d}_h(\Delta(c_{\mu,\lambda} \otimes c_{\mu',\lambda'}))
     & = \sum r_h(c_{\mu,\lambda_i} \otimes c_{\mu',\lambda'_j})
\otimes l_h(c_{-\lambda_i,\lambda} \otimes c_{-
\lambda_j',\lambda'})\\
     & = \sum\lambda_i(h)\lambda_j'(h)(-\lambda_i)(h)(-
\lambda_j')(h) c_{\mu,\lambda_i} \otimes c_{\mu',\lambda'_j}
\otimes c_{-\lambda_i,\lambda} \otimes c_{-\lambda_j',\lambda'}\\
     & =\Delta(c_{\mu,\lambda} \otimes c_{\mu',\lambda'})
\end{align*}
     The second assertion follows from the first because $\theta
\otimes m$ is a $\Gamma$-equivariant map.
\end{pf}

\begin{thm}
     The set $\cal{S}=\{1 \otimes k_{-\lambda} \mid \lambda \in
2\bL^+\}$ is an Ore subset of $\xi(\cqdg)$. The map $\xi$ extends
to an algebra isomorphism 
$$ \xi:  \cqdg_{\cal{S}} \to (\cqg \otimes \uqg^{\cop})^\Gamma.
$$
\end{thm}

\begin{pf}
We first show that  $\cal{S} \subset \xi(\cqdg)$. Notice that for
$\lambda \in \bL^+$, 
$$ \xi( c^\lambda_{\nu,\lambda} \otimes 1) = \sum
c^\lambda_{\nu,\mu_i} \otimes
l^+(c^\lambda_{-\mu_i,\lambda}) = c^\lambda_{\nu,\lambda}
\otimes k_{-\lambda}.$$
Similarly,
$$\xi(1 \otimes c^\lambda_{\nu,w_0\lambda}) =
c^\lambda_{\nu,w_0\lambda} \otimes k_{w_0 \lambda}.$$
Hence if 
$$\Delta(c^\lambda_{-\lambda,\lambda}) =
     c^\lambda_{-\lambda,\mu_i} \otimes
c^\lambda_{-\mu_i,\lambda}$$
then
$$1 = \epsilon(c^\lambda_{-\lambda,\lambda}) =
     \sum S(c^\lambda_{-\lambda,\mu_i})
c^\lambda_{-\mu_i,\lambda}$$
Hence
$$\sum \xi(1\otimes c^{-w_0\lambda}_{\mu_i,-\lambda})\xi
(c^\lambda_{-\mu_i,\lambda}\otimes 1)
     =\sum (c^{-w_0\lambda}_{\mu_i,-\lambda}
c^\lambda_{-\mu_i,\lambda})        
          \otimes k_{-2\lambda} = 1\otimes k_{-2\lambda} $$
as required.

     To prove the result it now remains to show that 
\begin{equation} \label{star}
\forall r \in (\cqg \otimes \uqg^\cop)^\Gamma, \;\;  \exists s \in
\cal{S} \text{ such that }rs \in \xi(\cqdg).
\end{equation}
 Notice that for $\lambda \in \bL$ and $h \in \Gamma$,
$$ d_h(1\otimes k_\lambda) = \lambda(h) 1\otimes k_\lambda. $$
From this it follows that
$$(\cqg \otimes \uqg^\cop)^\Gamma =(\cqg \otimes U^0)^\Gamma (1
\otimes \uqg^\cop)^\Gamma. $$
Since $\cal{S}$ is contained in the units of $\uqg$, it suffices
to verify condition \ref{star} for the two factors separately.

     By \cite[9.2.2]{Jo} $\cqg \otimes U^0$ is spanned by
elements of the form
$$c^\lambda_{\nu,w_0\lambda} c^{\lambda'}_{\nu',\lambda'} \otimes
k_\mu$$
where $\nu,\nu',\mu \in \bL$ and $\lambda, \lambda' \in \bL^+$.
This element is invariant if and only if $w_0\lambda
+\lambda'+\mu \in 2\bL$. On the other hand the image of $\cqdg$
contains 
$$\xi(1\otimes c^\lambda_{\nu,w_0\lambda}) \xi(
c^{\lambda'}_{\nu',\lambda'}\otimes 1) =
c^\lambda_{\nu,w_0\lambda} c^{\lambda'}_{\nu',\lambda'} \otimes
k_{w_0\lambda-\lambda'}.$$
Now for any $\gamma \in \bL$, there exists a $\gamma' \in \bL^+$
such that $\gamma'-\gamma \in \bL^+$. Thus 
if $w_0\lambda + \lambda' +\mu =2 \gamma$, then 
$$\mu - 2(\gamma'-w_0\lambda) = w_0\lambda -\lambda'
-2(\gamma'-\gamma).$$
Hence,
$$(c^\lambda_{\nu,w_0\lambda} c^{\lambda'}_{\nu',\lambda'}
\otimes k_\mu) 
(1 \otimes k_{-2(\gamma'-w_0\lambda)}) =
     \xi(c^\lambda_{\nu,w_0\lambda}\otimes
c^{\lambda'}_{\nu',\lambda'})
     (1\otimes k_{-2(\gamma'-\gamma)})
$$
which lies in $\xi(\cqdg)$.

     Notice that $(\uqg^{\cop})^\Gamma$ is generated over
$(U^0)^\Gamma$ by the elements $f_ik_i$ and $e_i$. Hence it
remains to verify condition \ref{star} for these elements. Let
$\lambda \in \bL^+$, let $\alpha_i$ be a simple root, set $\mu =
\lambda -\alpha_i$ and let $\nu \in \bL$. Then 
\begin{align*}
\xi(c^\lambda_{\nu,\mu}\otimes 1) & = \sum c^\lambda_{\nu,\mu_j}
\otimes l^+(c^\lambda_{-\mu_j,\mu})\\
     & =c^\lambda_{\nu,\lambda} \otimes
l^+(c^\lambda_{-\lambda,\mu}) + 
          c^\lambda_{\nu,\mu}\otimes l^+(c^\lambda_{-\mu,\mu}) 
\end{align*}
Now $l^+(c^\lambda_{-\lambda,\mu}) =f_ik_{-\mu}$ and
$c^\lambda_{\nu,\mu} \otimes l^+(c^\lambda_{-\mu,\mu})
\in (\cqg \otimes U^0)^\Gamma$. Hence by the paragraph above,
$\xi(\cqdg)$ contains $ c^\lambda_{\nu,\lambda} \otimes
f_ik_{-\mu}k_\gamma$ for some $\gamma \in -2\bL^+$. As noted
above,
$$\xi(1 \otimes c^{-w_0\lambda}_{\mu_j,-\lambda}) =
c^{-w_0\lambda}_{\mu_j,-\lambda} \otimes k_{-\lambda}.$$
Hence, for a suitably chosen $\gamma \in \bL^+$, $\xi(\cqdg)$
contains 
$$\sum (c^{-w_0\lambda}_{\mu_j,-\lambda} \otimes k_{-\lambda})
(c^\lambda_{-\mu_j,\lambda} \otimes f_ik_{-\mu}k_\gamma) =
( 1 \otimes f_ik_i)(1 \otimes k_{\gamma-2\lambda})$$
as required. A similar argument works for the element $e_i$.
\end{pf}
 
\section{Noetherianity}

     We now prove that the double of a standard quantum group
$\cqg$ is again Noetherian. We use the approach of Brown and
Goodearl \cite{BG}. Since both $\cqg$ and $\uqg$ are homomorphic
images of this algebra, this provides a straightforward and
unified proof that both $\cqg$ and $\uqg$ are noetherian. Recall
that if $\beta$ is a braiding on a Hopf algebra $A$, then the
braiding on the category of right comodules is given by the maps
$\beta_{V \otimes W} : V \otimes W \to W \otimes V$ where
$$\beta(v \otimes w) = \sum w_0 \otimes v_0 \beta(v_1, w_1)$$
(see, for instance, \cite{Ha}).

\begin{defn}
Let $(A, \beta)$ be a braided Hopf algebra and let $V$ be a
finite dimensional right
$A$-comodule. A full flag $0=V_0 \subset V_1 \subset \dots
\subset V_n=V$ of subspaces of $V$ is said to be 
{\em $\beta$-invariant} if
$$\beta_{V\otimes W}(V_i \otimes V) = V \otimes V_i$$
The flag is said to be {\em strongly $\beta$-invariant} if
$$\beta_{V\otimes W}(V_i \otimes W) = W \otimes V_i$$
for all right comodules $W$.
\end{defn}

     If $V$ is a right $A$-comodule with basis $v_i$ for $I \in
I$ and structure map $\rho:V \to V \otimes A$, then the subspace
of $A$ spanned by
$$\{ a_{ij} \mid \rho(v_j) = \sum v_i \otimes a_{ij} \}$$
is denoted by $C(V)$ \cite[p129]{Ab}.

\begin{defn} A right comodule $V$ over a Hopf algebra $A$ is said
to be a {\em generator} for $A$ if $A = k\langle C(V)\rangle$, or
equivalently, if $A = \sum_jC(V^{\otimes j})$.
\end{defn}

     One of the main results of \cite{BG} is the following
(slightly reworked into the language of braided Hopf algebras).

\begin{thm} Let $(A, \beta)$ be a braided Hopf algebra. Suppose
there exits a finite dimensional right $A$-comodule $V$ such that
\begin{enumerate}
\item $V$ is a generator for $A$;
\item $V$  has a $\beta $-invariant flag.
\end{enumerate}
Then $A$ is Noetherian.
\end{thm}

\begin{pf} See \cite[Theorem 4.4]{BG}
\end{pf}

     The key result is the following.

\begin{lem} All finite dimensional right $\cqdg$-comodules have
a strongly $\gamma'$-invariant flag.
\end{lem}

\begin{pf}
     It suffices to prove the result for comodules of the form $V
\otimes V'$ for comodules $V$ and $V'$ over $\cqg$. It follows
from the form of $\beta$ given above that there
exist strongly $\beta$-invariant flags of $V$ and $V'$ such that
$$ \beta_{V \otimes V'} (V_i \otimes V'_j ) = V'_j \otimes V_i$$
Now let $W$ and $W'$ be two other such $\cqg$-comodules. Notice
that 
$$\gamma'_{V\otimes V',W \otimes W'} =
     \beta_{14}(\tau \beta_{13})(\tau\beta^{-
1}_{42})(\tau\beta_{32}^{-1})$$ where $\beta_{kl}$ denotes the
map
induced from $\beta$ on the
$k$ and $l$-th components. Using this one can observe easily that
$$\gamma'_{V\otimes V',W \otimes W'}(V_i \otimes V'_j \otimes W
\otimes W')= W \otimes W'\otimes V_i \otimes V'_j $$
From this it follows easily that $V \otimes V'$ has a strongly
$\gamma'$-invariant flag.
\end{pf}

\begin{thm} The algebra $\cqdg$ is Noetherian for any connected
semi-simple algebraic group $G$.
\end{thm}

\begin{pf} It remains to notice that $\cqdg$ has a finite
dimensional right comodule
generator. The argument is the same as in the classical case.
\end{pf}

\begin{cor}
The algebras $\cqg$ and $\uqg$ are noetherian.
\end{cor}

\section{Further Non-standard quantum Groups}

We now discuss briefly  a generalization of double quantum
groups and a lifting theorem which provides a method of
constructing new families of  nonstandard quantum groups
corresponding to certain families of solutions of the modified
classical Yang-Baxter equation. 

We begin again with a very general result on lifting of Hopf
algebra twists. The proof is routine.

\begin{thm}\label{lift} Let $\phi: A \to B$ be a homomorphism of
Hopf algebras and let $\sigma$ be a 2-cocycle on $B$. Then
$\sigma'=\phi^*(\sigma)$ is a 2-cocycle. Moreover the induced map
$$\phi: A_{\sigma'} \to B_{\sigma}$$
between the corresponding twisted Hopf algebras is a Hopf algebra
homomorphism.

In particular, if $C$ is a braided Hopf  algebra and $\phi: A \to
C \otimes C$ is a Hopf algebra map, then there exists a 2-cocycle
$\sigma'$ on $A$ such that the map 
$$\phi: A_{\sigma'} \to C \Join C$$
is a Hopf  algebra homomorphism.
\end{thm}

In order to apply this result to standard quantum groups, one
needs
Hopf algebra maps of the form $\cqg \to \C_{q'}[G' \times G']$.
For
this it is not quite enough to find morphisms $G'\times G' \to
G$.
Let us say that a morphism of $G' \to G$ of connected,
simply-connected semi-simple
groups is {\em admissible} if the induced map on the Lie algebras
arises from a Dynkin diagram embedding as defined in
\cite[10.4.5]{Jo}. Then Braverman has shown that there is a Hopf
algebra surjection $\cqg \to \C_{q'}[G']$ and moreover that all
such maps arise in this way \cite{Br}.

\begin{thm}\label{qglift}
Let $G$ and $G'$ be connected, simply connected semi-simple
algebraic groups. Let $\psi:G'\times G'\to G$ be an admissible
embedding. Then there is a 2-cocycle $\sigma'$ on $\cqg$ and  an
associated non-standard quantum group
$\C_{q,\psi}[G]=\cqg_{\sigma'}$ such that
\begin{enumerate}
\item  $\C_{q,\psi}[G]$ is a braided Hopf algebra;
\item the category $\Comod$-$\C_{q,\psi}[G]$ is equivalent as a
braided rigid monoidal category to $\Comod$-$\cqg$.
\end{enumerate}
Moreover there is a natural surjective homomorphism,
$\C_{q,\psi}[G] \to \C_{q'}[D(G')]$ where $q'=q^r$ for some
rational number $r$ .
\end{thm}

For instance, if $2m\leq n$, then we have admissible
embeddings  $\psi: SL(m,\C) \times SL(m,\C) \to SL(n,\C)$. This
yields some nonstandard quantum groups of the form
$\C_{q,\psi}[SL(n)]$. Since the braiding $\gamma$ is essentially
the universal $T$-matrix, these quantizations appear to be
related to some of the ``esoteric quantum groups'' constructed by
Fronsdal and Galindo \cite{FG}.

The connection between these quantum groups and the solutions of
the modified classical Yang Baxter equation classified by Belavin
and Drinfeld appears to be the following. Let $\g$ and $\g'$ be
the
Lie algebras of $G$ and $G'$ respectively. Let $\bB$ and
$\bB'$ be bases for the roots of $g$ and $\g'$ respectively.
Then we may choose $\bB$ and $\bB'$ such that the map
$\psi$
induces an embedding of the Dynkin diagram $\bB'\times
\bB'$
into $\bB$. Let $\bB_1$ and $\bB_2$ be the images of
the
two copies of $\bB'$ in $\bB$ and let $\tau: \bB_1 \to
\bB_2$ be the natural isomorphism. Then $\tau$ is a triple in
the sense of \cite{BD}. We conjecture that $\C_{q,\psi}[G]$ can
be
regarded as a deformation in the algebraic sense \cite{DP} of the
algebra of functions on the group $G$ with Poisson structure
given
by a solution of the modified Yang-Baxter equation associated to
$\tau$. On the other hand, given any triple $\tau: \bB_1 \to
\bB_2$ in the sense of \cite{BD} which is disjoint in the
sense
that $(\alpha, \alpha')=0$ for all $\alpha \in \bB_1$ and
$\alpha' \in \bB_2$, then we may construct an admissible map
of
the form given in the theorem. Thus we may think of the quantum
group as being constructed directly from this data.

Notice that these cocycle twists are deceptively simple.
Although
the coalgebra structure is preserved by the twist, the algebra
structure is altered quite dramatically. Composing the
homomorphism
$\C_{q,\psi}[G] \to \C[D(G')]$ with the map $\C_{q'}[D(G')] \to
U_{q'}(\g')$ yields a Hopf algebra map $\C_{q,\psi}[G]\to
U_{q'}(\g')$. Thus $\C_{q,\psi}[G]$ has a category of finite
dimensional representations equivalent to those of $U_{q'}(\g')$.
This is quite different from the case of the standard quantum
groups where all finite dimensional representations are 
one-dimensional \cite[9.3.11]{Jo}. On the other hand the
existence of
this map is consistent with the philosophy on non-standard
quantum
groups advanced in \cite{H1}. Briefly the conjecture suggested by
this work is the following. Let $\C_q[G,r]$ be an algebraic
deformation of the group $G$ with Poisson structure given by a
solution $r$ of the modified classical Yang Baxter equation
associated to a triple $(\tau, \bB_1,\bB_2)$. Then there
should be a Hopf algebra homomorphism 
$\C_q[G,r] \to U_{q'}(\tilde{\g},\tilde{r})$ where $
U_{q'}(\tilde{\g},\tilde{r})$ is the FRT-dual of the quantization
$\C_{q'}[(\tilde{G},\tilde{r})$ of the reductive group
$\tilde{G}$
associated to the reductive Lie algebra $\tilde{\g}$ constructed
in
\cite[Theorem 6.4]{H1}. The semi-simple part of $\tilde{\g}$ is
the
semi-simple Lie algebra given by $\bB_1$.

Using Theorem \ref{lift} we may easily generalize the notion of
the
double $A\Join A$ of a braided Hopf algebra $A$ to a twisted
product $A^{\Join n}$ of $n$ copies of $A$ in the following way.
Using the map $m\otimes 1:(A\Join A)\otimes A \to A \otimes A$
and
Theorem \ref{lift} we may twist $(A\Join A)\otimes A$ to obtain a
new Hopf algebra $A^{\Join 3}$ which itself maps onto $A$.
Continuing in this way, we may iteratively construct $A^{\Join
(n+1)}$ from $A^{\Join n}$ using the map $A^{\Join n} \otimes A
\to
A \otimes A$ and Theorem \ref{lift}. Clearly we may also view the
construction $A^{\Join n}$ as a series of cocycle twists of the
$n$-fold tensor product $A^{\otimes n}= A \otimes A \otimes A
\dots
\otimes A$. Hence there exists a single cocycle $\tau$ such that
$A^{\Join n} \cong A^{\otimes n}_\tau$.

Again we may lift this construction to obtain interesting
families
of non-standard quantum groups. The algebra $\cqg^{\Join n}$ can
be
thought of as a nonstandard quantization of $G \times G \times
\dots \times G$. If $G$ and $G'$ are simply connected,  then any
admissible embedding $G' \times G' \times \dots \times G' \to G$
induces a homomorphism $ \cqg \to \ C_{q'}[G'\times \dots \times
G'] $
and
we may lift the twisting of $\C_q[G']^{\Join n}$ to a twisting on
$\cqg$ using Theorem \ref{lift}. Thus, for instance, the natural
block-diagonal embedding $\psi: SL(2,\C)\times SL(2,\C)\times
SL(2,\C)$ into $SL(6,\C)$ yields a new quantum group
$\C_{q,\psi}[SL(6)]$. 

One of the original motivations for this work was the desire to
find a construction of the Cremmer-Gervais quantum groups
\cite{H1}
from the standard quantum groups using a cocycle twist. However
none of the above constructions cover this case. The existence of
such a cocycle is equivalent to the existence of an equivalence
of braided rigid monoidal categories 
between the category of right comodules over the Cremmer-Gervais
quantum groups and the category of right comodules over the
standard quantum groups. Thus this remains a key problem in the
study of non-standard quantum groups. A positive answer to this
question would presumably also lead to a far more general
procedure
for constructing non-standard quantum groups.

\end{document}